\documentclass[12pt]{article}
\usepackage{amsmath,amssymb,amsthm}
\newtheorem{thm}{Theorem}
\newtheorem{prop}[thm]{Proposition}
\newtheorem{lem}[thm]{Lemma}
\newtheorem{cor}[thm]{Corollary}
\newtheorem{rem}[thm]{Remark}

\title{On the final limit of a transition matrix}
\author{Helmut Kahl}

\begin{document}
\maketitle	

\begin{abstract}
An efficient algorithm is presented for computation of the limit of exp(t A) for t towards infinity where A denotes an intensity matrix of finite dimension.
\end{abstract}

{\bf Mathematical Subject Cassification(2010).} 15B51, 15A16, 60J22, 60J27.


\section*{Introduction}
As a function of $t \in \mathbb{R}$ we call
\begin{equation*}
	t \mapsto \exp(t A) = \sum\limits_{n = 0}^{\infty} \frac{t^n}{n!} A^n
\end{equation*}
\textit{the fundamental matrix of generator} $A \in \mathbb{C}^{m \times m}$. Its limit $P$ for $t \to \infty$ is called the \textit{final limit} of $A$. The existence of $P$ depends on its eigenvalues only. A complex square matrix $A$ is called \textit{semistable} when its eigenvalues have negative real part or are zero, and if zero is an eigenvalue of $A$ it must be \textit{semisimple}, i.e. the dimension of its \textit{left nullspace} $N(A) := \lbrace x \in \mathbb{C}^m | x A = 0 \rbrace$ equals the algebraic multiplicity of eigenvalue zero. For describing $P$ the \textit{left image space} $R(A) := \lbrace x A \in \mathbb{C}^m | x \in \mathbb{C}^m \rbrace$ of $A$ is also important. In the following the zero matrix is denoted by $O$, a zero vector by $0$ and the identity matrix by $I$. We recall general facts about $P$, like given e.g. in \cite{Campbell-Rose}, Thm. 1.

\begin{thm}\label{thm_limit}
	The final limit $P$ of $A \in \mathbb{C}^{m \times m}$ exists if and only if $A$ is semistable. If $P$ exists $\mathbb{C}^m$ is the direct sum of $N(A)$ and $R(A)$, and it holds $P A = A P = O$. Hence $x \mapsto x P$ is the projector onto $N(A)$ along $R(A)$.
\end{thm}

By transposition an analogue holds for the right nullspace and image space. A \textit{row intensity matrix} $(a_{i j})_{i , j \in Z}$ is defined by $a_{i j} \ge 0$ for all $i , j \in Z$ with $i \ne j$ and by vanishing row sums $\sum\limits_{j \in Z} a_{i j} = 0$ for all $i \in Z$. Analogously one defines a \textit{column intensity matrix} by vanishing column sums. The fundamental matrix of an intensity matrix is called a \textit{transition matrix}. Namely, the probability distribution P of the Markov process $\lbrace X_t \rbrace_{t \ge 0}$ corresponding to a row intensity matrix $A$ on the \textit{state space} $Z$ with initial probability row vector $\pi$ (\textit{initial distribution}) fulfils
\begin{equation*}
	(\textnormal{P} \lbrace X_t = j \rbrace)_{j \in Z} = \pi \exp(t A) , t \ge 0 .
\end{equation*}
Hereby the off-diagonal entries of $A$ are the \textit{conditional intensities} of the Markov process. It is well known that, hereby, the transition matrix is \textit{row stochastic}, i.e. all its entries are nonnegative and all its row sums are equal to one. Analogously one defines a \textit{column stochastic} matrix. Hence we have

\begin{prop}\label{prop_bounded}
	The transition matrix is a function with constant row sum norm and thus bounded on the nonnegative real axis.
\end{prop}

In order to obtain interesting expectation values the main task is the computation of the \textit{final distribution}
\begin{equation*}
	\lim\limits_{t \to \infty} (\textnormal{P} \lbrace X_t = j \rbrace)_{j \in Z} = \pi \lim\limits_{t \to \infty}\exp(t A) = \pi P .
\end{equation*}
That $P$ exists is proven with probability theoretic means, e.g. in \cite{Stroock}, Thm. 5.4.6. That Theorem yields also a description of $P$ which conforms the description given in Theorem \ref{thm_compt}. Here we offer an existence proof relying on Proposition \ref{prop_bounded}.

\begin{cor}\label{cor_exist}
	The final limit of an intensity matrix $A$ exists and equals the projector onto $N(A)$ along $R(A)$.
\end{cor}
\begin{proof}
	According to \cite{Serre}, Prop. 5.12 all eigenvalues of $A = (a_{i j})$ are elements of the compact disc of radius $\rho := \max \lbrace -a_{i i} \mid i \in Z \rbrace$ around $-\rho$. Hence all eigenvalues of $A$ have negative real part or are zero. Assume $A x \ne 0$ and $A^2 x = 0$ for some column vector $x$. Then we have $A^n x = 0$ for all $n \ge 2$. It follows $\exp(t A) x = x + t A x$. But the latter function of $t$ is unbounded because of $A x \ne 0$. Hence also the transition matrix is unbounded which contradicts Proposition \ref{prop_bounded}. So for all column vectors $x$ the equation $A^2 x = 0$ implies $A x = 0$. So due to \cite{Serre}, Thm. 3.4 zero is semisimple. Now the assertion follows by Theorem \ref{thm_limit}.
\end{proof}

\section*{The final limit of a stochastic matrix}
According to \cite{Campbell-Rose}, Thm. 1 the final limit of an intensity matrix $A$ equals $I - A A^{\sharp}$ where $A^{\sharp}$ denotes the $\textit{Drazin-inverse}$ of $A$; for definition see e.g. \cite{Meyer}, eq. (5.10.6). There is an analogue about stochastic matrices $Q$: Due to \cite{Meyer}, ch. 8.4 the limit
\begin{equation*}
	Q_{\infty} := \lim\limits_{n \to \infty} \frac{1}{n} (Q + Q^2 + ... + Q^n)
\end{equation*}
exists and equals $I - A A^{\sharp}$ for $A := Q - I$. Hence we obtain:

\begin{cor}\label{cor_stoch}
	For a stochastic matrix $Q$ we have
\begin{equation*}
	Q_{\infty} = \lim\limits_{t \to \infty} \exp(t (Q - I)) .
\end{equation*}
\end{cor}

For a time-discrete Markov process $(X_n)_{n \in \mathbb{N}_0}$, also called a \textit{Markov chain}, it holds
\begin{equation*}
	(\textnormal{P} \lbrace X_n = j \rbrace)_{j \in Z} = \pi Q^n , n \in \mathbb{N}_0
\end{equation*}
where $\pi$ is again the initial distribution and $Q$ the constant transition matrix. Because of Corollary \ref{cor_exist} the problem of finding its final distribution $\pi Q_{\infty}$ can be regarded as the special case of finding the final distribution of a time-continuous Markov process with intensity matrix $Q - I$. So an efficient computation of $P$ becomes important also for Markov chains.

\section*{The final limit of an intensity matrix}
We show how to compute via simple matrix operations the final limit of a row intensity matrix. The other case is treated analogously with columns instead of rows or also can be reduced to the former case by transposition. The next fact establishes an isomorphism from the left nullspace of a row intensity matrix onto the left one-eigenspace of a certain row stochastic matrix $Q$. For the corresponding Markov process $Q$ is the \textit{transition matrix of the embedded Markov chain} $(Y_n)_{n \in \mathbb{N}_0}$ where $Y_n$ is the random state at the discrete time $n$. So stopping times are ignored as always in case of a Markov chain.

\begin{lem}\label{lem_embed}
	For a row intensity matrix $A = (a_{i j})$ let $D$ be the diagonal matrix with entries $d_i := -a_{i i}$. In case $d_i = 0$ substitute $d_i$ by one, so that $D$ becomes invertible. Then $Q := I + D^{-1} A$ is row stochastic, and the function $x \mapsto x D$ is an isomorphism from the left null-eigenspace of $A$ onto the left one-eigenspace of $Q$.
\end{lem}
\begin{proof}
	By Definition the entries of $Q = (q_{i j})_{i , j \in Z}$ are as follows: For $i \in Z$ we have $q_{i i} := 0$ in case $a_{i i} \ne 0$ and $q_{i i} := 1$ otherwise. For $i , j \in Z$ with $i \ne j$ it holds $q_{i j} := a_{i j} / d_i$. Thus $Q$ is row stochastic. For a row vector $x$ the equation $x B = 0$ is equivalent with $(x D) Q = x D$. This implies the latter assertion.
\end{proof}

According to \cite{Serre}, ch. 3.11 a square matrix $(a_{i j})$ over a commutative integral domain is called \textit{reducible} when there is some non-empty, proper subset $I$ of the index set $Z$ such that for all $i \in I$ and all $j \in Z \setminus I$ it holds $a_{i j} = 0$. Otherwise the matrix is called \textit{irreducible}. Note that $A$ is irreducible if and only if its transpose $A^{\textnormal{T}}$ is irreducible. The following fact can be seen as a version of the Perron-Frobenius-Theorem (see \cite{Serre}, Thm. 8.2).

\begin{prop}\label{prop_irred}
	An irreducible row intensity matrix $A$ has a one-dimensional nullspace and there is a unique row vector $p := (p_j)_{j \in Z}$ defined by $p A = 0$, $p_j > 0$ for all $j \in Z$ and
	\begin{equation*}
		\sum\limits_{j \in Z} p_j = 1 .
	\end{equation*}
\end{prop}
\begin{proof}
	The first assertion is shown in \cite{Serre}, ch. 8.3.3. The spectral radius of a stochastic matrix is one due to \cite{Serre}, ch. 8.5. Hence the second assertion follows from \cite{Serre}, Thm. 8.2 and Lemma \ref{lem_embed}.
\end{proof}

From now on we use the notation $A_{I, J} := (a_{i j})_{i \in I , j \in J}$ for nonempty subsets $I , J$ of $Z$ and $A_J := A_{J, J}$. An index $r \in Z$ is called \textit{recurrent with respect to} $A$ when there is some $r$ containing subset $J$ of $Z$ which for the submatrix $A_J$ is irreducible. Then $J$ is called the \textit{recurrence class} of $r$. A non-recurrent index is called \textit{transient}. Then for the set $R$ of recurrent indices $A_R$ has, after suitable permutation of $Z$, block diagonal form with blocks $A_J$. And for the set $T := Z \setminus R$ of transient indices we have $A_{R , T} = 0$ in case $T \ne \emptyset$.

\begin{rem}\label{rem_classes}
	
	a) The finiteness of $Z$ guarantees the existence of a recurrent index and so of a recurrence class with respect to any matrix $A = (a_{i j})_{i , j \in Z}$. In case $a_{i j} \ge 0$ for all $i , j \in Z$ let
	\begin{equation*}
		(c_{i j})_{i , j \in Z} := A + A^2 + ... + A^{n-1} .
	\end{equation*}
	Then index $i$ is in the same recurrence class as index $j \ne i$ if and only if $c_{i j}$ and $c_{j i}$ do not vanish. This explains how to compute recurrence classes and the set of transient states also for an intensity matrix, since its diagonal elements may be ignored hereby.
	
	b) Due to Proposition \ref{prop_irred} for a recurrence class $J \subset Z$ with respect to a row intensity matrix $A$ there is a unique row vector $p_J := (p_j)_{j \in Z}$ defined by
	\begin{equation*}
		p_J A = 0, p_j > 0 \textnormal{ for } j \in J, p_j = 0  \textnormal{ for } j \in Z \setminus J  \textnormal{ and } \sum\limits_{j \in Z} p_j = 1 .
	\end{equation*}
	The vectors $p_J$ form a basis of the nullspace of $A$. In case $T = \emptyset$ this is already clear by the first assertion of Proposition \ref{prop_irred} and by the block diagonal form mentioned above. And in case $T \ne \emptyset$ it follows by the following Proposition.
\end{rem}

\begin{prop}\label{prop_invert}
	If the set $T$ of transient indices with respect to an intensity matrix $A$ is not empty then  $A_T$ is invertible.
\end{prop}
\begin{proof}
	The definition of transience implies the existence of some $i \in T$ and some $j \in Z \setminus T$ with $a_{i j} \ne 0$. By transposition we may assume without loss of generality that $A$ is a row intensity matrix. The entry $q_{i j}$ of the row stochastic matrix $Q$ in Lemma \ref{lem_embed} fulfils also $q_{i j} \ne 0$. Hence $Q_T$ is not row stochastic but entriewise smaller or equal to a row stochastic matrix. Therefore $Q_T$ does not have eigenvalue one according to \cite{Serre}, ch. 8.3.2, Lemma 12. So by the isomorphism of Lemma \ref{lem_embed} the left null space of $A_T$ is zero.
\end{proof}

The following theorem gives a constructive description of the final limit which exists according to Corollary \ref{cor_exist}.

\begin{thm}\label{thm_compt}
	For a row intensity matrix $A$ with finite index set $Z$ let $P$ be the final limit of its transition matrix. For every recurrence class $J \subset Z$ with respect to $A$ a row of $P$ with index in $J$ equals the vector $p_J$ in Remark \ref{rem_classes}b). In case the set $T \subset Z$ of transient indices with respect to $A$ is not empty the matrices $P_T$ and $A_{T , R} P_R + A_T P_{T , R}$ with $R := Z \setminus T$ are zero. Hereby $A_T$ is invertible, whence
	\begin{equation*}
		P_{T , R} = -A_T^{-1} A_{T , R} P_R .
	\end{equation*}
	For the column vector
	\begin{equation*}
		(f_{i , J})_{i \in T} := - A_T^{-1} \sum\limits_{j \in J} (a_{i j})_{i \in T}
	\end{equation*}
	 and the set $\mathfrak{C}$ of recurrence classes the row of index $i \in T$ of $P$ equals
	 \begin{equation*}
	 	\sum\limits_{J \in \mathfrak{C}} f_{i , J} p_J .
	 \end{equation*}
\end{thm}
\begin{proof}
	Let us define a matrix $P$ such that its row of index $i \in J$ equals $p_J$ and, if $T \ne \emptyset$, such that $P_T = O$ and $A_{T , Z} P = A_{T , R} P_R + A_T P_{T , R} = O$. The latter equation is possible according to Proposition \ref{prop_invert}. From the definition of a recurrence class $J$ it follows $a_{i j} = 0$ for $i \in J$ and $j \in Z \setminus J$. Hence we have $A_{J , Z \setminus J} P_{Z \setminus J , Z} = O$. Since all row vectors $p_J$ of $P$ of index $i \in J$ are equal and the row sums of $A$ vanish we have $A_{R , J} P_{J , Z} = A_{J , J} P_{J , Z} = O$. It follows $A_{R , Z} P = O$. Together with $A_{T , Z} P = O$ by definition of $P$ we obtain $A P = O$. Let $x$ be a row vector with $x A = 0$. According to Remark \ref{rem_classes}b) $x$ is a linear combination of the $p_J$. Because of $p_J P = p_J$ for all $J \in \mathfrak{C}$ it follows $x P = x$. Thus $P$ is the projector onto $N(A)$ along $R(A)$. So $P$ is the final limit according to Corollary \ref{cor_exist}. For $j \in J \in \mathfrak{C}$ it holds
	\begin{equation*}
		A_T P_{T , \lbrace j \rbrace} = - A_{T , R} P_{R , \lbrace j \rbrace} = - p_{j j} \sum\limits_{j \in J} A_{T , \lbrace j \rbrace} = p_{j j} A_T (f_{i , J})_{i \in T} .
	\end{equation*}
	By cancelling out $A_T$ we obtain $p_{i j} = f_{i , J} p_{j j}$ for $i \in T$ and hence the identity for the $i$-th row of $P$.\footnote{$f_{i , J}$ is the enter probability from the transient state $i$ into the recurrence class $J$.}
\end{proof}

\section*{Computation of the final limit}
Essentially, the computational costs of the following pseudo-algorithm suggested by Theorem \ref{thm_compt} are an exponentiation of the $n \times n$ matrix $A$ to the power of $n-1$ like explained in Remark \ref{rem_classes}a) and (in the worst case of irreducible $A$) solving a regular linear equation system of dimension $n$. Thus it may be called efficient.

\begin{rem}\label{rem_compt}
	Computation of the final limit $P$ from row intensity matrix $A$ with index set $Z$:
	\begin{itemize}
		\item[0.] Initialise $P =(p_{i j}) := O$ (zero matrix of dimension like $A$) and $T := Z$.
		\item[1.] Compute the set $\mathfrak{C}$ of recurrence classes with respect to $A$.
		\item[2.] For $J \in \mathfrak{C}$:
		\begin{itemize}
			\item[2.1.] Compute $p = (p_j)_{j \in J}$ with $p A_J = 0$ and $\sum\limits_{j \in J} p_j = 1$.
			\item[2.2.] Substitute every row of $P_J$ by $p$.
			\item[2.3.] Substitute $T$ by $T \setminus J$.
		\end{itemize}
		\item[3.] If $T$ is not empty:
		\begin{itemize}
			\item[3.1.] For $J \in \mathfrak{C}$:
			\begin{itemize}
				\item[3.1.1.] Compute $f = (f_i)_{i \in T}$ s.t. $A_T f + A_{T , J} (1 ... 1)^{\textnormal{T}} = 0$.
				\item[3.1.2.] For $i \in T$: For $j \in J$: Substitute $p_{i j}$ by $f_i p_{j j}$.
			\end{itemize}
		\end{itemize}
		\item[4.] Output $P$.
	\end{itemize}
\end{rem}

\end{document}